\documentstyle[12pt]{article}
\title{On the asymptotic geometry of the hyperbolic plane}
\author {Iosif Polterovich\thanks{Department of Theoretical Mathematics,
The Weizmann Institute of Science, Rehovot 76100, Israel; 
Department of Mechanics and
Mathematics, Moscow State University, Moscow 119899, Russia; e-mail:
iossif@wisdom.weizmann.ac.il.}\and
Alexander Shnirelman\thanks{School of 
Mathematical Sciences, The Raymond and Beverly Sackler
Faculty of Exact Sciences, Tel Aviv University, Ramat Aviv 69978,
Israel; e-mail: 
shnirelm@math.tau.ac.il. This author is supported by 
the French-Israeli Scientific and Technical
Cooperation Program PICS No. 296.}}
\begin{document}
\maketitle
\abstract {Asymptotic subcone of an unbounded metric space is another
metric space, capturing the structure of the original space at
infinity. In this paper we define a functional metric space $S$
which is an asymptotic subcone of the hyperbolic plane. 
This space is a real tree branching at every its point. Moreover,
it is a homogeneous metric space such that any real tree with
countably many vertices can be isometrically embedded into it.
This implies that every such tree is also an asymptotic subcone of the 
hyperbolic plane.} 

\def \phi{\varphi}
\def \epsilon{\varepsilon}
\newtheorem{definition}{Definition}[section]
\newtheorem{lemma}{Lemma}[section]
\newtheorem{theorem}{Theorem}[section]
\newtheorem{corollary}{Corollary}[section]
\section{Introduction and main results}
 
Let $(X, d_X)$ be a metric space $X$ with a distance function $d_X$. 
Suppose that $X$ has infinite diameter, i.e. the function $d_X$ is unbounded.
Then one may ask what is the structure of the space $X$ "at infinity".
Intuitively, structure at infinity is what is seen if one looks at
the space $X$ from an infinitely far point (see [Gr2]).
M. Gromov suggested several ways to treat this notion rigorously.
In this paper we follow one of them.

\proclaim Definition 1.1.
{\rm Let $(X, d_X)$ be a metric space with an infinite diameter.
A metric space $(T, d_T)$ can be} isometrically embedded at infinity {\rm into
the space $X$ if for every point $t \in T$ there exists an
infinite sequence $\{x_t^i\}$, $i=1,2,..$ of points in $X$,
such that for some  fixed sequence of positive $\epsilon_i \to 0$
$$\lim_{i\to\infty}\epsilon_i\cdot d_X(x_{t_1}^i, x_{t_2}^i)=
d_T(t_1, t_2) 
\eqno (1.2)
$$  
for every $t_1,t_2 \in T$.}

In other words, to every point of the space
$T$ we may put in correspondence a sequence of points in $X$ going to infinity,
such that the ``normalized'' pairwise distances between the sequences in $X$
tend to the distances between the corresponding points in $T$. 

The Definition 1.1. somehow clarifies what the ``structure at infinity''
means but it is still too difficult to work with it.

In order to proceed we need to define a certain 
class of ``geometrically simpler'' metric
spaces --- {\it geodesic} metric spaces (see [Gr1], [GhH]):

\proclaim Definition 1.3.
{\rm Let $x_0,x_1$ be two points of the metric space $X$ and let  
$a=d_X(x_0,x_1)$ be the distance between them. A} geodesic segment {\rm in
$X$ {\rm connecting
$x_0$ and $x_1$} is an isometric inclusion $g:\lbrack 0,a \rbrack \to X$ such that
$g(0)=x_0$, $g(a)=x_1$. The image of this inclusion sometimes is also 
called a geodesic segment.
A metric space $X$ is} geodesic {\rm if for every two points
$x_1,x_2 \in X$ there exists a 
(not necessarily unique) geodesic segment connecting these two points.}   

For example, every complete Riemannian manifold is a geodesic space
due to Hopf-Rinow theorem. All metric spaces which appear in this paper 
are geodesic.

\proclaim Definition 1.4.
{\rm A metric space $X_0$ is an} asymptotic subcone {\rm of the space $X$
if $X_0$ is geodesic and its every finite subset of points can be
isometrically embedded at infinity into the space $X$.}
 
Asymptotic subcones were introduced in [Gr1] in the context of 
{\it hyperbolic}  metric spaces (see [Gr1], [GhH]). There are several
equivalent formal definitions of a hyperbolic metric space but all of them
demand additional non--trivial geometric explanations.  
We are not presenting any of them 
them here since throughout this paper we deal with the familiar 
Lobachveskian (hyperbolic) plane which is  the most well--known
example of a hyperbolic metric space. 

In fact, there is another very well--known class of hyperbolic
spaces~--- these are $0$-hyperbolic spaces (see [Gr1], [GhH] for the
definition of $\delta$-hyperbolicity) , or {\it real
trees}:

\proclaim Definition 1.5.
{\rm A metric space $X$  is a} real tree {\rm if it satisfies the following
conditions:
\noindent 1) For every two distinct 
points of the space there exists a unique geodesic segment joining them.
\noindent 2)If two geodesic segments
$[a,b], [b,c]$ have exactly one endpoint  $b$ in
common,
their union is also a geodesic segment.}
 
Real trees play the key role in the asymptotic geometry of
the hyperbolic metric spaces.
There is a general theorem that every asymptotic
subcone of a hyperbolic space is a real tree (see [GhH]). 
Moreover, as it was stated by Gromov (see [Gr1], [GhH]), 
if every asymptotic subcone of a metric space is a real tree
than this space is hyperbolic. Therefore, one may define hyperbolic
spaces as metric spaces whose all asymptotic subcones are real trees.
For the hyperbolic groups Gromov considered this definition as the 
most intuitive ([Gr1]).

The property of being a real tree does not give a complete description
of a metric space. Different real trees can be very much unlike each
other. Therefore in order to describe asymptotic subcones of a particular
hyperbolic metric space (a Lobachevskian plane in our case)
it is not satisfactory to say they are just real trees,  a  far more 
``explicit'' construction is desirable.

We found  only one example of an asymptotic subcone of a hyperbolic
plane in the literature --- a star--shaped tree formed by
$k$ segments with a common vertex (see [GhH]). Clearly, the asymptotic
geometry of a hyperbolic plane is much richer.
Hyperbolic plane is a homogeneous metric space. It is quite natural
to look for asymptotic subcones sharing this property. Such a subcone
should be real tree branching at every its point --- already an object 
which is quite
difficult to imagine. We also want our subcone to contain
all ``simple'' examples of asymtotic subcones (it follows from the 
Definition 1.4. that every geodesic subset of an asymptotic subcone 
of the space $X$ is itself an asymptotic subcone of $X$). As a criterion
of ``simplicity'' we choose countability of the number of vertices of
the real tree:
 
\proclaim Definition 1.6.
{\rm A real tree is } thick {\rm if it allows an isometric
inclusion of any real tree with countably many vertices.}

Thus, we want to ``materialize'' an asymptotic subcone of a hyperbolic
plane which is a homogeneous and thick real tree. The surprising fact is
that such a substantial part of the ``structure at infinity'' of the
Lobachevskian plane can be described as a certain simple functional
space.

\proclaim Definition 1.7.
{\rm Let $S$ be 
the set of all continuous real functions $f(t)$ 
defined on a finite interval 
$\lbrack 0,\rho \rbrack$, $0 \le \rho< \infty$ (each function has its
own $\rho$), such that 
$f(0)=0$ for all $f \in S$.
Define the following metric on $S$: 
$$d_S(f_1,f_2)=(\rho_1-s)+(\rho_2-s),
\eqno(1.8)
$$   
where 
$[0,\rho_i]$ is the domain of  the function $f_i$ , $i=1,2$ and
$$s=\sup \lbrace t|f_1(t')=f_2(t') 
\quad \forall t' < t 
\rbrace.
$$ 
This defines the} metric space S.
{\rm The number $s$ is called the} moment of segregation  
{\rm of the functions  $f_1(t), f_2(t)$.}

Let us formulate the main result of our paper:

\proclaim Theorem 1.9.
The space $S$ is a thick real tree and a homogeneous metric
space. It is an asymptotic subcone of the hyperbolic plane.
 
We prove this theorem in the next two sections. In the last section
we show that the completion of an asymptotic subcone of a metric
space is itself an asymptotic subcone of a metric space
(motivated by the fact that the space $S$ is non--complete). 
The paper is completed by two appendices.

\medskip

\noindent {\bf Remark 1.10.}
Most of the results of the present paper were announced in 
[PSh]. Some constructions introduced here were used in [Sh] to describe
the {\it asymptotic cone}, or the {\it asymptotic space}  (see [Gr0], [Gr2]) 
of the Lobachevskian plane by means of non--standard analysis 
(see [D]).  

\medskip

\noindent{\it Acknowledgements.}
The authors are grateful to M.~Gromov, L.~Poltero-
\noindent vich, V.~Buchstaber and A.~Vershik
for helpful discussions and support.

\bigskip

\bigskip

\bigskip

\bigskip

\section{Properties of the metric space $S$}
In this section we study some properties of the metric space $S$.

\proclaim Lemma 2.1.
The function $d_S$ defined by {\rm ( 1.8)} is a metric. 

\noindent {\bf Proof.}
We need to check that the metric (1.8) satisfies the triangle inequality.
Let $f_1,f_2,f_3$ be three functions in $S$, 
$\rho_1, \rho_2, \rho_3$ be the lenghts of their domains and
$s_{12}, s_{13},s_{23}$ be their segregation moments. 
We may always assume that $s_{12} \le s_{13} \le s_{23}$
Then clearly $s_{12}=s_{13}$. Therefore,
$$d_S(f_1,f_2)+d_S(f_1,f_3)=2\rho_1+\rho_2+\rho_3-4s_{12} \ge
\rho_2+\rho_3-2s_{12} \ge 
$$
$$
\ge \rho_2+\rho_3-2s_{23} = d_S(f_2,f_3).$$
The two other inequalities are proved similarly.

\medskip

\proclaim Lemma 2.2.
The space $S$ is a real tree.

\noindent {\bf Proof.}
Let $f_1, f_2 \in S$ be two arbitrary functions, $\rho_1, \rho_2$ be their
domains and $s$ be their moment of segregation. By (1.8)
$d_S(f_1,f_2)=\rho_1+\rho_2-2s$. 
Consider the following inclusion $g:\lbrack 0,\rho_1+\rho_2-2s \rbrack \to X$:
$$g(x)=\cases{
\lbrace f_1(t),\, 0\le t\le \rho_1-x \rbrace, 0\le x\le \rho_1-s; \cr
\lbrace f_2(t),\, 0\le t \le x+2s-\rho_1 \rbrace, \rho_1-s\le x \le 
\rho_1+\rho_2-2s.\cr}
 \eqno (2.3)
$$
This inclusion is isometric and clearly unique with such property, therefore
the condition 1) of the Definition 1.5. is verified. In order to check the
condition 2) we note that any two geodesic segments $[f,g]$, $[g,h]$ may have 
exactly one point $g$ in common if and only if the function $h$ 
is the extension of the function $g$ and segregates from it not earlier 
than from the function $f$,
or, symmetrically,
if the function $f$  is the extension of the function $g$ and segregates 
from it  not earlier than 
from the function $h$. In both cases the formula (2.3) implies that  
$\lbrack f,h \rbrack$ is also a geodesic segment. 
Therefore $S$ is a real tree which completes the proof.   

\medskip

\proclaim Lemma 2.4.
The space  $S$ is a thick tree.

\noindent {\bf Proof.} 
Let $T$ be an arbitrary real tree with countably many vertices. We ``brush''
this tree in the following way. Fix some isomorphism between natural
numbers and the set of all vertices. Let $a_{ij}$ be the distance between
vertices corresponding to the numbers $i$ and $j$, and $\lbrace k_n \rbrace$ 
be an infinite strictly increasing sequence of natural numbers.
Now we build a mapping from $T$ into $S$. Let the vertex $1$ go to zero
(by zero we denote the function defined 
and equal to 0 at the single point $0$).
The vertex $2$ goes to a linear function $f(t)=k_1t$ defined on the
interval $\lbrack 0, a_{12} \rbrack$. In order to find the image of the
vertex $3$ we find from $a_{12}$, $a_{13}$ and 
$a_{23}$ where it branches from $1$ and $2$;
let $s_3$ be the abscissa of this point. Therefore on the interval
$\lbrack 0, s_3 \rbrack$ it is already defined and on the interval
$\lbrack s_3, a_{13} \rbrack$ we set it  to be linear 
with the angular coefficient
$k_2$ (the free term is found from continuity). 
Repeating the same inductive algorithm for all $n$ (if $n-1$ vertices
are already built we find the abscissa  $s_n$ of its point of segregation
from the already built tree and continue the function by setting it linear with
the coefficient $k_{n-1}$ on the interval  
$\lbrack s_n, a_{1n} \rbrack$) we get an inclusion of $T$ into $S$.
It is isometric by construction since the sequence $\lbrace k_n \rbrace$ 
is strictly increasing  and hence segregation is defined correctly.

\medskip

Now let us prove that the space $S$ is homogeneous, i.e. 
for every two its points
there exists a one-to-one isometry moving one point to another.

\proclaim Lemma 2.5.
The metric space $S$ is homogeneous. 

\noindent {\bf Proof.}
Clearly it is sufficient to construct a one-to-one isometry $F$ which moves
any function  to zero. 
Denote the preimage of zero by $f_0(t)$, let $[0,\rho]$ be 
its domain. Let $f(t)$ be  any other function with the domain $[0,a+b]$, where
$a$ is the moment of segregation of the functions $f_0(t), f(t)$.
If $a < \rho$ then the image of $f$ is given by  the function 
$F(f(t))$,such that $F(f(t))=0$  on $\lbrack 0,\rho-a\rbrack$ and  
$F(f(t))=f(t-\rho+2a)-f_0(a)$ on $\lbrack \rho-a,\rho-a+b \rbrack$.
If $a=\rho$, i.e. the function $f$ is a ``continuation'' of $f_0$, the
construction is more complicated. Let us choose some infinite sequence
of continuous functions $\lbrace g_n(t) \rbrace$ such that $g_1(t)$
is identically zero and for any two elements of this sequence
their moment of segregation is zero. For example we can take the sequence
$$g_n(t)=\frac{(2^n-1)t}{2^n}, n=0,1,2,....$$
Consider the function 
$F^*(f(t))=f(t+\rho)-f_0(\rho)$
defined on $\lbrack 0,b \rbrack$. If there exists $0<c\le b$ such 
that $F^*(f(t))=g_n(t)$ identically on $\lbrack 0,c \rbrack$ 
for some $n\ge 0$ then $F(f(t))=F^*(f(t))+g_{n+1}(t)$ on
$\lbrack 0,b \rbrack$, otherwise simply $F(f(t))=F^*(f(t))$.
One can easily check that $F$ is indeed an isometry.

\medskip

As it was mentioned in the introduction, homogeneity is a very
important feature of the space $S$ since we are interested in it
as in the model of the asymptotic space of the hyperbolic plane,
and such a model can not be considered ``good'' if it does not preserve
such a fundamental property of the initial space. 
At the same time the question about the relations between the isometries
of $S$ and of the initial hyperbolic plane remains open.
 
\section{Asymptotic subcones of the hyperbolic plane}  

Let $X$ be the hyperbolic plane. For the conveniency of computations we use 
its Poincare unit disc model  
$$X=\lbrace x \in C : |x|<1 \rbrace.$$ 

For every point $x\in X$, denote by $\rho$ the non-Euclidean
distance between $x$ and $0$, the centre of $X$, and by $\phi$
the polar angle; thus, $(\rho,\phi)$ are the non-Euclidean
polar coordinates of the point $x$.

The Euclidean distance between $0$ and $x$ is denoted by $r$.

There is the following relation between $\rho$
and $r$ (see [Be]):
$$r=\frac{e^{\rho}-1}{e^{\rho}+1}.
\eqno(3.1)
$$

The distance between the two points $x_1,x_2 \in X$ is given by the
formula ([Be]):
$$d_X(x_1,x_2)=\ln \frac{|1-x_1\overline{x_2}|+|x_1-x_2|}
{|1-x_1\overline{x_2}|-|x_1-x_2|}  \eqno (3.2)$$
We rewrite this formula in Euclidean polar coordinates. If $x_1=(r_1,\phi_1)$,
$x_2=(r_2,\phi_2)$ then the formula (3.2) transforms into the
following:
$$d_X(x_1,x_2)=\ln \frac{\sqrt{\frac{1+(r_1r_2)^2-2r_1r_2\cos(\phi_1-\phi_2)}
{r_1^2+r_2^2-2r_1r_2\cos(\phi_1-\phi_2)}}+1}
{\sqrt{\frac{1+(r_1r_2)^2-2r_1r_2\cos(\phi_1-\phi_2)}
{r_1^2+r_2^2-2r_1r_2\cos(\phi_1-\phi_2)}}-1}
\eqno (3.3)
$$

After a number of elementary transformations
of (3.3) we get:
$$d_X(x_1,x_2)=\ln \frac{1+A}{1-A}, \eqno(3.4)$$
where
$$A^2=1-\frac{8}{(2-\beta^2)(t+1/t)^2+\beta^2(s+1/s)^2},$$
$$
\beta^2=1-\cos(\phi_1-\phi_2),\, s^2=e^{\rho_1+\rho_2},\,t^2=e^{\rho_1-\rho_2}.
$$  

In order to prove that the space $S$ is an asymptotic subcone of
the hyperbolic plane we introduce an auxilary metric space $D$:

\proclaim Definition 3.5.
{\rm Consider
the set of all real functions $f(t)$ 
defined on a finite semi-interval 
$\lbrack 0,\rho)$, $0 \le \rho< \infty$ ($\rho$ also depends on the function),
such that $f(0)=0$ and $f(t)=0$ everywhere but at a finite number of points.
Denote the }metric space $D$ {\rm by endowing this set with the metric
$d_D$, whose expression is given by the formula (1.8).}

The space $D$ is also a real tree (the proof is exactly the same as 
of the Lemma 2.2.). In fact, it can be isometrically included into the
space $S$ but is not isometric it (and hence ``smaller'' than the space $S$).
This statement is proved in the Appendix A.

\proclaim Lemma 3.6.
The space $D$ is an asymptotic subcone of the hyperbolic plane.

\noindent{\bf Proof.}
Let $f_1(t), f_2(t)$ be two arbitrary  functions from $D$, let 
$\rho_1 \ge \rho_2$ be
the lengths of their domains, $\lbrace s_{i,1}, ..., s_{i,N_i}\rbrace$
be their supports and $f_i(s_{i,k})=a_{i,k}$, $k=1,..N_i$,  $i=1,2$.  
Consider the following two sequences of points of the hyperbolic plane:
$x_i(n)=(\rho_i(n), \phi_i(n))$, where $\rho_i(n)=\rho_i/\epsilon_n$,  
$\phi_i(n)=\sum_{k=0}^{N_i} a_{i,k}e^{-s_{i,k}/\epsilon_n}$, $i=1,2$,
where $\lbrace \epsilon_n \rbrace$ is an arbitrary sequence of positive 
numbers tending to zero.
Then after a simple asymptotic analysis of the formula (3.4) we get 
 
$$
\lim_{n \to \infty} \epsilon_n d_X(x_1(n),x_2(n))= \lim_{n \to \infty}
\epsilon_n\ln(e^{(\rho_1-\rho_2)/\epsilon_n}+
e^{-2s/\epsilon_n}e^{(\rho_1+\rho_2)\epsilon_n})=
$$
$$
=\rho_1+\rho_2-2s,
$$
where $s$ is exactly the moment of segregation of the functions 
$f_1(t), f_2(t)$. Therefore,
$$\lim_{n \to \infty}\epsilon_n d_X(x_1(n),x_2(n))=d_D(f_1,f_2).$$
Comparing this
with (1.2) and recalling that 
$D$ is a geodesic space being a real tree completes the proof of our lemma.

\medskip
We call $D$ the {\it discrete subcone} of the hyperbolic plane.

Now we are able to complete the proof of our main result.

\proclaim Theorem 3.7.
The space $S$ is an asymptotic subcone of the hyperbolic
plane.

\noindent{\bf Proof.}
Let $f_k(t)$, $k=1,...,n$ be an arbitrary finite subset of the
space $S$. Denote by $[0,\rho_k]$ the domains of the functions
$f_k(t)$, $k=1,..n$ and by   
$s_{k_1k_2}$ the moments of segregation of the
functions $f_{k_1}$ and $f_{k_2}$,  $k_1,k_2=1,..n$.
We choose $\frac{n(n-1)}{2}$ pairwise distinct 
points $t^1_{k_1k_2}$ on the positive half-line
such that
$$ f_{k_1}(t_{k_1k_2}) \not= f_{k_2}(t_{k_1k_2})
\quad {\rm and} \quad t^1_{k_1k_2}-s_{k_1k_2}< 1/4 
\eqno(3.8)
$$ 
for every two different $k_1, k_2 = 1..n$ (the first condition should
be checked only if $t^1_{k_1k_2} \le \min(\rho_{k_1},\rho_{k_2})$
since otherwise it does not make sense).
Let us put down these points in the growing order and denote the resulting 
ordered set
by $\lbrace t^1_j \rbrace$: $t^1_{j_1}> t^1_{j_2}$ for all $j_1>j_2$,
$j_1,j_2=1,..\frac{n(n-1)}{2}$.  
Take an arbitrary sequence 
of positive numbers 
$\lbrace \epsilon_i \rbrace$, $i=1,2,..$
tending  to zero. 
Consider the following sequences of points of the hyperbolic plane:

$$x^1_{i,k}=(\rho_k/\varepsilon_i, \phi^1_{i,k}(\epsilon_i)),$$
where
$$ \phi^1_{i,k}(\epsilon_i)= \sum_{j=1}^{J^1_k} e^{-t^1_j/\varepsilon_i}
f_k(t^1_j), \qquad k=1,..,n ,$$
where $J^1_k$  is the number of points $t^1_j$ which lie in the domain
$ [0,\rho_k]$ of the function $f_k$; if there are no such points set 
$\phi^1_{i,k}(\epsilon_i)\equiv 0$.  

Denote by $g^1_k(t)$, $k=1,...n$, the functions 
belonging to the discrete subcone
$D$ such that each $g^1_k(t)$ has the same domain $[0,\rho_k]$ as $f_k(t)$ and
$$
g^1_k(t)=\cases{f_k(t),\quad  t=t^1_j, j=1,...J^1_k;\cr
0,\quad  {\rm otherwise}.\cr}
$$
 
Then, by lemma 3.6 there exists a number 
$I_1$, such that for all $k_1,k_2=1,..n$:
$$
|\varepsilon_{I_1} d_X (x^1_{I_1,k_1},x^1_{I_1,k_2})-
d_D(g^1_{k_1},g^1_{k_2})| \le 1/2.
\eqno (3.9)
$$ 
Therefore, by the choice of the points $t^1_{k_1k_2}$ we have:
$$
|\varepsilon_{I_1} d_X (x^1_{I_1,k_1},x^1_{I_1,k_2})-
d_S(f_{k_1},f_{k_2})| \le 
$$
$$
\le |\varepsilon_{I_1} d_X (x^1_{I_1,k_1},x^1_{I_1,k_2})-
d_D(g^1_{k_1},g^1_{k_2})|+|d_D(g^1_{k_1},g^1_{k_2})-d_S(f_{k_1},f_{k_2})|\le
$$
$$
\le 1/2+ 2|t^1_{k_1k_2}-s_{k_1k_2}|<1/2+2\cdot 1/4 = 1.
$$ 

Now we take new points $t^2_{k_1k_2}$ which satisfy the condition (3.8)
with $1/8$ instead of $1/4$ in the right--hand side and get the new set 
of points $\lbrace t^2_j \rbrace$, the new
sequences $x^2_{i,k}$ and the new  functions
$g^2_k(t) \in D$, $k=1,..n$.

Similarly, there exists a number $I_2$, greater than $I_1$, such that
$$
|\varepsilon_{I_2} d_X (x^2_{I_2,k_1},x^2_{I_2,k_2})-
d_D(g^2_{k_1},g^2_{k_2})| \le 1/4
$$
for all $k_1,k_2=1,..,n$, and therefore
$$
|\varepsilon_{I_2} d_X (x^2_{I_2,k_1},x^2_{I_2,k_2})-
d_S(f_{k_1},f_{k_2})| \le 1/4+2|t^2_{k_1k_2}-s_{k_1k_2}|<1/2
$$
for all $k_1,k_2=1,..,n$.

Continuing this procedure (which is possible due to the Lemma 3.6. and
the assumption that $f_k(t)$ are
continuous functions)
we get the sequence $\lbrace I_N \rbrace$, 
$I_N \to \infty$ as $N \to \infty$, 
the sequences $\lbrace {\varepsilon_{I_N} \rbrace \to 0}$
(or just  $\lbrace \varepsilon_N \rbrace$),
$\lbrace x^N_{I_N,k} \rbrace$  (or just $\lbrace x_{N,k} \rbrace$),
and the functions $g^N_k(t) \in D$, $k=1,..,n$
such that 
$$
|\varepsilon_N d_X (x_{N,k_1},x_{N,k_2})-
d_S(f_{k_1},f_{k_2})| \le 
$$
$$
\le
|\varepsilon_N d_X (x_{N,k_1},x_{N,k_2})-
d_D(g^N_{k_1},g^N_{k_2})|+
2|t^N_{k_1k_2}-s_{k_1k_2}|<
$$
$$
< 1/2^N+2 \cdot 1/2^{N+1}=1/2^{N-1}
$$

Thus we have proved that there exists a limit
$$
\lim_{N \to \infty} \varepsilon_N d_X(x_{N,k_1},x_{N,k_2})=
d_S(f_{k_1},f_{k_2})
$$
for all $k_1,k_2= 1,..,n.$
Therefore $S$ is indeed an asymptotic subcone of the 
hyperbolic  plane.
  
\smallskip

The Theorem 3.7. together with the Lemmas 2.2., 2.4. and 2.5.
completes the proof of the Theorem 1.9.

We call $S$ the {\it continuous subcone} of the hyperbolic plane.
In fact it is  an asymptotic subcone of an
infinitely narrow neighborhood of a single half-line on the hyperbolic plane,
as follows from the construction of the discrete subcone $D$.

Some modifications of the construction of the
continous subcone $S$ are considered in the Appendix B.

\medskip

We would like to conclude this section with a simple corollary of
the main result  which however reflects the ``wealth''
of the space $S$:

\proclaim Corollary 3.10.
Any real tree with countably many vertices is an asymptotic subcone of 
the hyperbolic plane.

\noindent{\bf Proof.}
Indeed, every geodesic subspace of an asymptotic subcone 
is itself an asymptotic subcone. Real trees are geodesic spaces
and by Lemma 2.4 any
real tree with countably many vertices can be isometrically inluded into 
the continuous subcone $S$.
Therefore every such tree is an asymptotic subcone of the hyperbolic plane.
  
\section{Completions of the asymptotic subcones}

The metric spaces $S$ and $D$ are non--complete           
metric spaces. For $D$ this is obvious;
for $S$ it follows from the following example.
Consider a sequence $\{f_k(t)\}$ of functions, defined on the segment
$[0,1-2^{-k}]$, and equal to  $f_k(t)=\sin(1/(1-t))$ on its domain of
definition. Clearly, this sequence has no limit 
in the space $S$ as $k \to \infty$.

Therefore, it is reasonable to consider the completion 
$\bar S$ 
of the  
continuous subcone $S$. 
Though it does not have such a simple functional description as the
original one, it is  also an asymptotic subcone of the hyperbolic plane
due to the following  simple general theorem (since we could not find it in
the literature we found it appropriate to state it here):

\proclaim Theorem 4.1.
A completion of an asymptotic subcone of a metric space is also
an asymptotic subcone of this metric space.

\noindent {\bf Proof.}
Let $X$ be our metric space, $X_0$ --- its asymptotic subcone and $\bar X_0$ 
--- the completion of $X_0$.
Let $f_1,f_2,...,f_k$ be a finite number of points in $\bar X_0$.
By definition, 
$$f_i=\lim_{n \to \infty} \phi_i^n,$$ 
where $\phi_i^n \in X_0$, $n \in {\bf N}$, $i=1,...,k$.
 
Choose a number $N_1=N_1(1/4)$ such that 
$$d_{\bar X_0}(\phi_i^n,f_i)<1/4
$$ 
for all $i=1,..,k$, $n \ge N_1$.
Then for all $l,m=1,..,k$ and $n \ge N_1$ we get

$$|d_{\bar X_0}(f_l,f_m)-d_{X_0}(\phi_l^n,\phi_m^n)|\le
d_{\bar X_0}(\phi_l^n,f_l)+d_{\bar X_0}(\phi_m^n,f_m) <
$$
$$
< \frac{1}{4}+\frac{1}{4}=
\frac{1}{2}
\eqno(4.2)
$$ 

The space $X_0$ is an asymptotic subcone of $X$. Therefore for its
finite subset $\phi_i^{N_1}$, $i=1,..,k$ there exists a 
sequence $\{ \epsilon_j \}$ of positive numbers tending to zero,
and $k$ sequences of points of the space $X$,  
$\{ x_j^{i,N_1} \}$, $i=1,..,k$,  such that for some $J_1=J_1(N_1,1/2)$
$$ |\epsilon_j d_X(x_j^{l,N_1},x_j^{m,N_1}) - 
d_{X_0}(\phi_l^{N_1},\phi_m^{N_1})| < \frac{1}{2}
\qquad \forall j\ge J_1.
\eqno (4.3)
$$

Therefore, by (4.2) and (4.3) we get:
$$ |\epsilon_{J_1} d_X(x_{J_1}^{l,N_1},x_{J_1}^{m,N_1}) - 
d_{\bar X_0}(f_l, f_m)| < \frac{1}{2} + \frac{1}{2}=1
$$
for all $l,m=1,..,k$.

Next, we choose  $N_2=N_2(1/8)$ and $J_2=J_2(N_2, 1/4)$
and similarly obtain
$$ |\epsilon_{J_2} d_X(x_{J_2}^{l,N_2},x_{J_2}^{m,N_2}) - 
d_{\bar X_0}(f_l, f_m)| < \frac{1}{4} + \frac{1}{4}=\frac{1}{2}
$$
for all $l,m=1,..,k$.
Let us note that we can always choose $N_2$ and $J_2$ in such a way that
$N_2 > N_1$ and $\epsilon_{J_2}< \epsilon_{J_1}/2$. 

Continuing this procedure analogously we get 
the sequence $\{\epsilon_{J_r}\}$ of positive numbers tending to zero 
and $k$ sequences of points in $X$, 
$\{x_{J_r}^{i,N_r}\}$, $i=1,..,k$, such that

$$ |\epsilon_{J_r} d_X(x_{J_r}^{l,N_r},x_{J_r}^{m,N_r}) - 
d_{\bar X_0}(f_l, f_m)| < \frac{1}{2^r} + \frac{1}{2^r}=\frac{1}{2^{r-1}}
$$
for all $l,m=1,..,k$,  which implies

$$\lim_{r \to \infty}\epsilon_{J_r} d_X(x_{J_r}^{l,N_r},x_{J_r}^{m,N_r})=
d_{\bar X_0}(f_l, f_m)|, \quad l,m=1,..,k.
\eqno(4.4)
$$
 
The relation (4.4) exactly means that $\bar X_0$ is an asymptotic subcone
of the metric space $X$, which completes the proof of the theorem. 
  
\bigskip

\section*{Appendix A}

We have shown that the spaces $S$ and $D$
are real trees. These
trees are branching at every point and the cardinal number of 
the set of their vertices is continuum (for $D$ this is
clear and for $S$ it follows from the fact that every 
continuous function is defined by its values at the rational points). 
As a metric space $S$ is ``larger'' than $D$, as it is shown by the following

\proclaim Lemma A.1.
The space $D$ can be isometrically included into the space 
$S$ but is not isometric to it.

\noindent {\bf Proof.} 
Let us show that $D$ can be 
isometrically included into $S$.
Let $\gamma(t)$ be any element of $D$, 
$Z_\gamma=\lbrace a_1<a_2<..<a_N \rbrace$ be the set  
of  points where $\gamma (t)$ is non-zero, 
and $[0,\rho_0)$ be the domain of $\gamma(t)$.
Set $a_0=0$, $a_{N+1}=\rho_0$.
Consider the following map $F: D\to S$: 
$$
F[\gamma(t)](t)=\sum_{i=0}^k \gamma(a_i)(t-a_i),\, \, a_k \le t \le a_{k+1},\, 
k=0,..,N.
\eqno (A.2)
$$
It is easy to see that $F$ is an isometric inclusion of $D$ into $S$.

Now, assume that $\Phi$ is a an isometric inclusion of
$S$ into $D$ and let $\Phi[f(t)](t)=\gamma(t)$
for some $f(t)\in S$  where $\gamma(t)\in D$ 
is defined as above.

Take some $0<\epsilon_0 < \rho_0-m_0$, where $m_0=a_N$, and consider two
functions $g_1(t)$, $h_1(t)$ 
in the $\epsilon_0$-neighborhood of $f(t)$ in $S$, such that none
of the functions $f(t)$, $g_1(t)$, $h_1(t)$
is the extension of any of the other two. Denote
$\alpha_1(t)=\Phi[g_1(t)](t)$, $\beta_1(t)=\Phi[h_1(t)](t)$. 

Let us show that 
$Z_\gamma \subset Z_{\alpha_1}$, $Z_\gamma \subset Z_{\beta_1}$.
If, for instance, $Z_\gamma$ is not a subset of $Z_{\alpha_1}$
then $d_D(\gamma,\alpha_1)\ge \rho_0-m_0>\epsilon_0$, 
which contradicts with the choice of $\alpha_1(t)$; the similar argument
is valid for $\beta_1(t)$.  
At least one of these inclusions is proper; 
indeed, if $Z_{\alpha_1}=Z_{\beta_1}= Z_\gamma$ then one 
of the functions $\alpha_1(t)$,
$\beta_1(t)$, $\gamma(t)$ is the  extension of the two others which is
impossible since $\Phi$ is an isometry and none of the functions
$f(t)$, $g_1(t)$, $h_1(t)$
is the extension of any of the other two.
Let this proper inclusion be $Z_\gamma \subseteq Z_{\alpha_1}$ .
Then the set $Z_{\alpha_1}$ consists of at least $N+1$ points.

Repeat this construction taking $g_1$ instead of $f$, and instead of  
$\epsilon_0$ taking  
$0< \epsilon_1< \min(\epsilon_0/2,\rho_1-m_1)$, where
$[0,\rho_1)$ is the domain of $\alpha_1(t)$ and 
$m_1$ is the maximal element in  $Z_{\alpha_1}$.
Similarly, we shall get the new function $\alpha_2(t) \in D$ such that
$Z_{\alpha_2}$ consists of at least $N+2$ points.

Continuing this procedure we obtain a sequence $\lbrace g_i(t)\rbrace$ in 
$S$ and the corresponding sequence $ \lbrace \alpha_i(t)\rbrace$ in $D$ 
such that
$Z_{\alpha_i}$ consists of at least $N+i$ points.
When $i \to \infty$ the functions $g_i(t) \to g(t)$, where $g(t) \in S $ since 
$\epsilon_i<\epsilon_0/2^i$ for all $i=1,2,..$ and therefore
the length of  the domain of $g(t)$ is finite --- it is not greater than
$r+2\epsilon_0$, where $r$ is the length of the domain of $f(t)$. 
Consider the image of $g(t)$ under the isometry $\Phi$;
it is equal to $\alpha(t)=\lim_{i\to \infty}\alpha_i(t)$.
By our construction $Z_\alpha \supset Z_{\alpha_i}$ for all $i=1,2,..$ ,
therefore $Z_\alpha$ consists of an infinite number of points. 
But this contradicts with the fact that $\alpha(t) \in D$,
and hence $S$ and $D$ are non-isometric.
This completes the proof of our lemma.      
 
\section*{Appendix B}

Instead of continuous functions with bounded
domain one could take generalized functions of bounded domain with the
distance defined by (1.8).
Such generalized functions are of finite order (see [GeS]), i.e. they
can be represented as finite sums of generalized derivatives of continuous
functions. One may check that integration preserving the condition $f(0)=0$
is an isometry with respect to our metric
(compare this with the formula (A.2.) --- 
in fact we have represented each element
of the discrete subcone $D$ as a finite sum of $\delta$-functions and 
intergrated twice). Hence every space of all generalized
functions of order less than some finite $m > 0$ is isometric to $S$:
its isometric inclusion into $S$ is obtained  by integrating $m$ times
every its element as described above, 
and isometric inclusion of $S$ into such 
space can be given by  an identical map.   
Similarly, if we take $C^m$-smooth functions,
we also get an isometric space. Therefore, all such functional spaces are
isometric asymptotic subcones of the hyperbolic plane. 

The construction of $S$ can be also generalized for
the Lobache\noindent vskian  space of arbitrary dimension $n$. 
In this case instead of  scalar functions one has to take 
$n$-component vector functions.

\bigskip

\bigskip

\bigskip

\bigskip

\medskip

\centerline{ {\large\bf References}}

\medskip

\noindent [Be] A. Beardon, The geometry of discrete groups, Springer--Verlag, 
1983.

\smallskip

\noindent  [D] M.~Davis, Applied Nonstandard Analysis, 
J. Wiley \& Sons, New York, 1977

\smallskip

\noindent [GeS] I.M. Gelfand, G.E. Shilov, 
Generalized functions, Vol. 2, Academic Press, 1968. 

\smallskip 

\noindent [GhH] E. Ghys, P. de la Harpe, Sur le groupes hyperboliques apres
Mikhael Gromov, Birkh\"auser, 1990. 

\smallskip

\noindent [Gr0] M. Gromov, Groups of polynomial growth and expanding
maps, IHES Math. Publ., N 53, 1981, 53-71.

\smallskip
 
\noindent [Gr1] M. Gromov, Hyperbolic Groups, in: Essays in group theory,
ed. S.M.Gersten, M.S.R.I. Publ. 8 , Springer--Verlag, 1987, 75-263.

\noindent [Gr2] M. Gromov, Asymptotic invariants of infinite groups, 
Geometric group theory. Vol. 2 (Sussex, 1991), London Math. Soc. Lecture
Note Ser. 182, Cambridge Univ. Press, 1993, 1-295.

\smallskip
\noindent [PSh] I. Polterovich, A. Shnirelman, An asymptotic subcone       
of the Lobachevskii plane as a space of functions, Russian Math. Surveys,  
v. 52 No. 4, 1997, 842-843.

\noindent [Sh]  A. Shnirelman, On the structure of  asymptotic space       
of the Loba-

\noindent chevsky plane, 1-23, to appear in the Amer. J. Math.

\end{document}